\documentclass[a4paper,12pt]{amsart}

\usepackage{amssymb} 
\usepackage{latexsym} 
\usepackage{amsfonts} 
\usepackage{amsmath} 
\usepackage{eucal} 
\usepackage{bm} 
\usepackage{bbm} 
\usepackage{graphicx} 
\usepackage[english]{varioref} 
\usepackage[nice]{nicefrac} 
\usepackage[all]{xy}

\newcommand{\N}{\mathbbm{N}}

\newcommand{\Hom}{{\rm Hom}}

\newcommand{\Ho}[1]{\textrm{\textnormal{H}}^{0}(#1)}

\usepackage{amsthm}

\theoremstyle{plain}
\newtheorem{thm}{Theorem}[section] 
\newtheorem*{thm*}{Theorem} 
 \newtheorem{prop}[thm]{Proposition}
 \newtheorem{lemma}[thm]{Lemma}
 \newtheorem{cor}[thm]{Corollary}
 \newtheorem{remark}[thm]{Remark}

\theoremstyle{definition}

\newtheorem{example}[thm]{Example}

\theoremstyle{remark}

\newtheorem*{remark*}{Remark}
\newtheorem*{claim*}{Claim}

\begin{document}

\title[On compactifications of the Steinberg zero-fiber]
{On compactifications of the Steinberg zero-fiber}
\author{Thomas Haahr Lynderup \and Jesper Funch Thomsen}
\address{Institut for matematiske fag\\Aarhus Universitet\\ \AA rhus, Denmark}
\email{haahr@imf.au.dk, funch@imf.au.dk}

\begin{abstract}
Let $G$ be a connected semisimple linear algebraic group over an
algebraically closed field $k$ of positive characteristic and let
$X$ denote an equivariant embedding of $G$. We define 
a distinguished Steinberg fiber $N$ in $G$, called the zero-fiber, 
and prove that the closure of $N$ within $X$  
is normal and Cohen-Macaulay. Furthermore, when $X$
is smooth we prove that the closure of $N$ is a 
local complete  intersection. 
\end{abstract}

\maketitle

\section{Introduction}
Let $G$ be a connected semisimple linear algebraic group over an
algebraically closed field $k$ of positive characteristic. The set 
of elements in $G$ with semisimple part within a fixed $G$-conjugacy 
class is called a Steinberg fiber. Examples of Steinberg fibers 
include the unipotent variety and the conjugacy class of a regular 
semisimple element. Lately there has been some interest in describing 
the closure of Steinberg fibers within equivariant embeddings of 
the group $G$ (see \cite{He}, \cite{HT}, \cite{Spr}, \cite{T}). 

In this paper we study the closure of a distinguished Steinberg 
fiber $N$ called the Steinberg zero-fiber (see Section 
\ref{Steinbergfibers} for the precise definition of $N$).
We will prove that the closure $\bar{N}$ of $N$ within any 
equivariant embedding $X$ of $G$ will be normal and 
Cohen-Macaulay. Moreover, when $X$ is smooth we will 
prove that $\bar{N}$ is a local complete intersection. 
These results will all be proved by Frobenius splitting 
techniques. As a byproduct we find that $\bar{N}$ has
a canonical Frobenius splitting and hence the set of global 
sections of any $G$-linearized line bundle on $\bar{N}$
will admit a good filtration.

The presentation in this paper is close to \cite{T}, 
but the setup is somehow opposite. More precisely, 
in loc.cit. the group $G$ was fixed to be of simply 
connected type while the Steinberg
fiber was arbitrary; in the present paper the Steinberg 
fiber is fixed but the semisimple group $G$ is of arbitrary 
type. 

It is worth noticing (see \cite{T}) that for a fixed 
equivariant embedding $X$ of $G$ the boundary 
$\overline{F}-F$ of the closure of a Steinberg fiber $F$ in $G$ 
is independent of $F$. This shows that the 
results obtained in this paper for the 
rather special Steinberg zero-fiber will
provide some knowledge about closures
of Steinberg fibers in general.
E.g. the boundary $\overline{F}-F$ 
will always be Frobenius split (see Theorem \ref{thmNbarFsplit}).
This suggests, that the results in this paper may
be generalized to arbitrary Steinberg fibers. 
However, we give an example (see Example \ref{example})
showing that this is not always the case.

\section{Notation}
Let $G$ denote a connected semisimple linear algebraic 
group over an algebraically closed field $k$. The
associated groups of simply connected and adjoint 
type will be denoted by $G_{\rm sc}$ and $G_{\rm ad}$
respectively. Let $T$ denote a maximal torus in $G$
and let $B$ denote a Borel subgroup of $G$ containing
$T$. The associated maximal torus and Borel subgroup 
in $G_{\rm ad}$ (resp. $G_{\rm sc}$) will be denoted
by $T_{\rm ad}$ and $B_{\rm ad}$ (resp. $T_{\rm sc}$
and $B_{\rm sc}$).

The set of roots associated to $T$ is denoted by $R$. 
We define the set of positive roots $R^+$ to be the nonzero 
$T$-weights of the Lie algebra of $B$. The set of 
simple roots $(\alpha_i)_{i \in I}$ will be indexed 
by $I$ and will have cardinality $l$. To each simple 
root $\alpha_i$ we let $s_i$ denote the associated 
simple reflection in the Weyl group $W$ defined by
$T$. The length of an element $w$ in $W$ is defined
as the number of simple reflections in a reduced 
expression of $w$. The unique element in $W$ of 
maximal length will be denoted by $w_0$. 

The weight lattice $\Lambda(R)$ of the root
system associated to $G$ is identified with the 
character group $X^*(T_{\rm sc})$ of $T_{\rm sc}$ 
and contains the set of $T$-characters $X^*(T)$. We
let $\alpha_i^\vee$, $i \in I$, be the set of simple coroots and 
let $\langle, \rangle$ denote the pairing between
coweights and weights in the root system $R$. A 
weight $\lambda$ is then dominant if $\langle \lambda
, \alpha_i^\vee \rangle$ for all $i \in I$. The 
fundamental dominant weight associated to $\alpha_i$,
$i \in I$, will be denote by $\omega_i$.

For a dominant weight $\lambda \in \Lambda(R)$ we 
let ${\rm H}(\lambda)$ denote the dual Weyl 
$G_{\rm sc}$-module with heighest weight 
$\lambda$, i.e. containing a $B$-semiinvariant 
element $v_\lambda^+$ of weight $\lambda$. 
The Picard group of $\nicefrac{G}{B}$ 
may be identified with the weight lattice 
$\Lambda(R)$ and we let $\mathcal{L}(\lambda)$ 
denote the line bundle associated to $\lambda
\in \Lambda(R)$. The line bundle 
$\mathcal{L}(\lambda)$ has a unique 
$G_{\rm sc}$-linearization so we may regard
the set of global sections of $\mathcal{L}(\lambda)$ 
as a $G_{\rm sc}$-module. We assume that the
notation is chosen such that the set of global 
sections of $\mathcal{L}(\lambda)$ is 
isomorphic to ${\rm H}(-w_0\lambda)$.
For $\lambda, \mu \in  \Lambda(R)$ we 
denote by $\mathcal{L}(\lambda,\mu)$
the line bundle $\mathcal{L}(\lambda)
\boxtimes \mathcal{L}(\mu)$ on the 
variety $\nicefrac{G}{B} \times
\nicefrac{G}{B}$.

\section{Steinberg fibers}
\label{Steinbergfibers}

The set of elements $g$ in $G$ with semisimple 
part $g_s$ in a fixed $G$-conjugacy class is called a 
{\it Steinberg fiber}. Any Steinberg fiber is a 
closed irreducible subset of $G$ of codimension
$l$ (see \cite[Thm.6.11]{St}). Examples of Steinberg
fibers include the conjugacy classes of regular 
semisimple elements and the unipotent variety;
i.e. the set of elements with $g_s$ equal to 
the identity element $e$. 

When $G=G_{\rm sc}$ is 
simply connected the Steinberg fibers may also 
be described as genuine fibers of
a morphism $ G_{\rm sc} \rightarrow k^l$. 
Here the $i$-th coordinate map is given by 
the $G_{\rm sc}$-character of the representation 
${\rm H}(\omega_i)$. In this formulation the 
fiber $N_{\rm sc}$ above $(0,0,\dots, 0)$ 
is called the {\it Steinberg zero-fiber}. 
When $G$ is arbitrary we define the Steinberg 
zero-fiber $N$ of $G$ to be the image of 
$N_{\rm sc}$ under the natural morphism 
$\pi : G_{\rm sc} \rightarrow G$.

The structure of the Steinberg zero-fiber is very 
dependent on the characteristic of $k$. In most cases 
$N$ is just the conjugacy class of a regular 
semisimple element. At the other extreme we
can have that $N$ coincides with the unipotent
variety of $G$.    

  \begin{remark}
    In the following cases the Steinberg zero-fiber and the
    unipotent variety of $G$ coincide :

     Type $\mathsf{A}_{n}$: when $n=p^{m}-1$ ($m \in \N$) and
     $p= {\rm char}(k) > 0$.

    Type $\mathsf{C}_{n}$: when $n=2^{m}-1$ ($m \in \N$) and
    ${\rm char}(k)=2$.

    Type $\mathsf{D}_{n}$: when $n=2^{m}$ ($m \in \N$) and
    ${\rm char}(k)=2$.

    Type $\mathsf{E}_{6}$: when ${\rm char}(k)=3$.

    Type $\mathsf{E}_{8}$: when ${\rm char}(k)=31$.

    Type $\mathsf{F}_{4}$: when ${\rm char}(k)=13$.

    Type $\mathsf{G}_{2}$: when ${\rm char}(k)=7$.
  \end{remark}

\section{Equivariant Embeddings}
An \emph{equivariant embedding} of $G$ is a normal $G \times G$-variety
containing a $G \times G$-invariant open dense subset isomorphic to 
$G$ in such a way that the induced $G \times G$-action on $G$ is 
by left and right translation.

\subsection{The wonderful compactification}
\label{wonderful}

When $G = G_{\rm ad}$ is of adjoint type there exists a distinguished 
equivariant embedding $\bm X$ of $G$ called the wonderful compactification
(see e.g. \cite{DCS}). 
The wonderful compactification of $G$ is a smooth projective 
variety with finitely many $G \times G$-obits $O_J$
indexed by the subsets $J$ of the simple roots $I$. We 
let ${\bm X}_{J}$ denote the closure of $O_J$ and assume
that the index set is chosen such that 
${\bm X}_{J'} \cap {\bm X}_{J} =   {\bm X}_{J' \cap J}$
for all $J,J' \subset I$. Then ${\bm Y} := {\bm X}_{\emptyset}$ 
is the unique closed orbit in ${\bm X}$. It is known that ${\bm Y}$
is isomorphic to $\nicefrac{G}{B} \times \nicefrac{G}{B}$ 
as a $G \times G$-variety.  

To each dominant element $\lambda$ in the weight lattice $\Lambda(R)$
we let 
$$ \rho_\lambda^{\rm ad} : G \rightarrow {\mathbb P}({\rm End}({\rm H}
(\lambda))),$$
denote the $G \times G$-equivariant morphism defined by 
letting $\rho_\lambda^{\rm ad}(e)$ be the element in  
${\mathbb P}({\rm End}({\rm H}(\lambda)))$ represented 
by the identity map on ${\rm H}(\lambda)$. Then it is
known (see \cite{DCS}) that $\rho_\lambda^{\rm ad}$ extends to a
morphism ${\bm X} \rightarrow {\mathbb P}({\rm End}({\rm H}
(\lambda)))$ which we also denote by $\rho_\lambda^{\rm ad}$. 

\begin{lemma}
\label{closedorbit}
Let $v_\lambda^+$ (resp. $u_\lambda^+$) denote a nonzero 
$B$-stable element in ${\rm H}(\lambda)$ (resp. ${\rm H}
(\lambda)^*$) of weight $\lambda$ (resp. $-w_0 \lambda$).
Identify ${\rm End}({\rm H}(\lambda))$ with ${\rm H}(\lambda) 
\boxtimes {\rm H}(\lambda)^*$. Then the restriction of 
$\rho_\lambda^{\rm ad}$ to $\bm Y \simeq \nicefrac{G}{B} 
\times \nicefrac{G}{B}$ is given by 
$$\rho_\lambda^{\rm ad} ( gB, g'B) = (g,g')[v_\lambda^+
\boxtimes u_\lambda^+].$$
\end{lemma}
\begin{proof}
It suffices to prove that the only $B \times B$-invariant 
element of the image of $\rho_\lambda^{\rm ad}$ is
$[v_\lambda^+\boxtimes u_\lambda^+]$. So let 
$x = [f]$ denote a $B \times B$-invariant
element of the image of $\rho_\lambda^{\rm ad}$
represented by an element $ f \in {\rm End}({\rm H}
(\lambda))$. Then $f$ is  $B \times B$-semiinvariant 
and thus when writing $f$ as an element of 
${\rm H}(\lambda) \boxtimes   {\rm H}(\lambda)^*$
it will be equal to $v_\lambda^+ \boxtimes v$ for some 
$B$-semiinvariant element $v$ of ${\rm H}(\lambda)^*$.

Let ${L}$ denote the unique simple submodule 
of  ${\rm H}(\lambda)$ and let $M$ denote 
the kernel of the associated morphism
${\rm H}(\lambda)^* \rightarrow L^*$. Assume 
that $v$ is contained in $M$. Then every 
$G \times G$-translate of $f$ is contained
in the subset $L \boxtimes M$ consisting
of nilpotent endomorphisms. Now we apply 
\cite[Lemma.6.1.4]{BK}. It follows that 
the closure $C$ of the $T \times T$-orbit 
through $\rho_\lambda^{\rm ad}(e)$ 
will contain an element represented by
a nilpotent endomorphism. But clearly 
every element in $C$ will be represented
by a semisimple endomorphism of 
${\rm H}(\lambda)$. This is a 
contradiction. 

As a consequence $v$ is not contained in $M$ 
and therefore its image in $L^*$ will be a 
nonzero $B$-semiinvariant vector. As a 
consequence $v$ must be a nonzero multiple
of $u_\lambda^+$.
\end{proof}

\subsection{Toroidal embeddings}
\label{sectiontoroidalembeddings} 
    
Now let $G$ be an arbitrary connected semisimple group. 
An equivariant embedding $X$ of $G$ is called toroidal 
if the natural map $\pi_{\rm ad} : G \rightarrow G_{\rm ad}$ extends to 
a morphism $X \rightarrow {\bm X}$. In the present paper 
toroidal embeddings will play a central role.  
 This is due to the following fact (see \cite[Prop.3]{Rit})

\begin{thm}
\label{thmprojectiveresolution}
Let $X$ be an arbitrary equivariant embedding of $G$. Then 
there exists a smooth toroidal embedding $X'$ of $G$ and a 
birational projective morphism $X' \rightarrow X$ extending 
the identity map on $G$. 
\end{thm}

An important property of toroidal embeddings is that 
for each dominant weight $\lambda$ there exists a 
$G \times G$-equivariant morphism
$$ \rho_\lambda :X \rightarrow {\mathbb P}({\rm End}({\rm H}
(\lambda)))$$ 
induced by $\rho_\lambda^{\rm ad}$. When $X$ is a 
complete toroidal embedding of $G$ and $Y$ is a closed 
$G \times G$-orbit in $X$ we may describe the restricted
mapping $Y \rightarrow {\mathbb P}({\rm End}({\rm H}
(\lambda)))$ using Lemma \ref{closedorbit}. 
Notice that  $Y$  maps surjectively to $\bm Y
\simeq \nicefrac{G}{B} \times \nicefrac{G}{B}$ and
that $Y$ is a quotient of $G \times G$. Consequently
$Y$ maps bijectively onto $\bm Y$ and hence $Y$ must 
be $G \times G$-equivariantly isomorphic to
$\nicefrac{G}{B} \times \nicefrac{G}{B}$. Moreover  

\begin{lemma}
\label{closedorbit2}
Let $X$ be a complete toroidal embedding of a connected
semisimple group $G$ and let $Y$ be a closed $G \times G$-orbit
in $X$. Then $Y$ is $G \times G$-equivariantly isomorphic 
to $\nicefrac{G}{B} \times \nicefrac{G}{B}$ and 
$$\rho_{\lambda}( gB, g'B)  \mapsto  (g,g')[v_\lambda^+
\boxtimes u_\lambda^+].$$
Consequently, the pull back of the ample generator 
$\mathcal{O}_\lambda(1)$ of the Picard group
of ${\mathbb P}({\rm End}({\rm H}(\lambda)))$ to 
$Y$ is isomorphic to $\mathcal{L}(\lambda, 
-w_0 \lambda)$.
\end{lemma}

\subsection{The dualizing sheaf of equivariant embeddings}

Let $X$ be a smooth equivariant embedding of $G$. As $G$
is an affine variety the complement $X - G$ of $G$ is of 
pure codimension 1 in $X$. Let $X_1, \dots,
X_n$ denote the irreducible components of $X - G$ which
are then divisors in $X$. Let $D_i$, $i=1,\dots,l$, denote
the closures of the Bruhat cells $B w_0 s_i B$ within
$X$. Then also $D_i$ is a divisor in $X$. When $X$ is 
the wonderful compactification $\bm X$ of
a group of adjoint type we will also denote $X_j$ 
and $D_i$ by ${\bm X}_j$ and ${\bm D}_i$ respectively.

\begin{prop}\cite[Prop.6.2.6]{BK}
\label{canonical}
The canonical divisor of the smooth equivariant embedding 
$X$ is 
$$K_X = - 2 \sum_{i \in I} D_i - \sum_{j=1}^n X_j.$$  
\end{prop}

The line bundle $\mathcal{L}(D_i)$ associated to the 
divisor $D_i$ is connected to the morphisms 
$\rho_\lambda$ as explained by 

\begin{lemma}
\label{sections}
Assume that $X$ is a toroidal embedding. Then 
there exists an isomorphism 
$$\mathcal{L}(D_{i}) \simeq \rho_{\omega_i}^*( 
\mathcal{O}_{\omega_i}(1)),$$
such that $\rho_{\omega_i}^*(u_{\omega_i}^+\boxtimes 
v_{\omega_i}^+)$ is the canonical section of
$\mathcal{L}(D_{i})$.
\end{lemma}
\begin{proof}
Consider first the case when $X$ is the wonderful
compactification of $G_{\rm ad}$. Consider the 
pull back $s_{\rm ad} :=
(\rho_{\omega_i}^{\rm ad})^*(u_{\omega_i}^+\boxtimes v_{\omega_i}^+)$   
of the global section  $u_{\omega_i}^+\boxtimes v_{\omega_i}^+$ 
of $\mathcal{O}_{\omega_i}(1)$. Then the zero divisor 
$(s_{\rm ad})_0$ of $s_{\rm ad}$ is $B \times B$-invariant. 
Thus there exist nonnegative integers $a_r$ and $b_j$ , 
for $r,j=1,\dots,l$, such that
$$ (s_{\rm ad})_0 = \sum_{r=1}^l a_r {\bm D}_r + \sum_{j=1}^l
b_j {\bm X}_j.$$
If $b_j>0$ for some $j$ then $s_{\rm ad}$ vanishes 
on $\bm Y$ which by Lemma \ref{closedorbit} is a 
contradiction. Hence, $(s_{\rm ad})_0 = \sum_{r=1}^l a_r {\bm D}_r$.
It is known (see e.g. \cite[Prop. 6.1.11]{BK}) that the 
restriction of $\mathcal{L}({\bm D}_i)$ to $\bm Y$ is isomorphic to 
${\mathcal L}(\omega_i, - w_0 \omega_i)$, so using using 
Lemma \ref{closedorbit2} it follows that $a_{i} = 1$ and 
$a_r = 0$ for $r \neq {i}$. This proves the statement when
$X$ is the wonderful compactification of $G_{\rm ad}$.

Consider now an arbitrary toroidal equivariant embedding 
$X$ of $G$. Let $s := \rho_{\omega_i}^*(u_{\omega_i}^+\boxtimes 
v_{\omega_i}^+)$  be the pull back of the the global section 
$u_{\omega_i}^+\boxtimes v_{\omega_i}^+$ of $\mathcal{O}_{\omega_i}(1)$. 
Then the zero divisor $(s)_0$ of $s$ is $B \times B$-invariant
and by the already proved case above we conclude that 
$$(s)_0 = c D_{i} + \sum_{j=1}^n d_j X_j,$$
for certain nonnegative integers $c>0$ and $d_j$, 
$j =1, \dots, n$. Assume $d_j >0$ for some $j$. Then  
$s$ vanishes on a $G \times G$-stable subset and 
as $s$ is the pull back of $s_{\rm ad}$ from $\bm X$ to $X$,
we conclude that $s_{\rm ad}$ also vanishes on a 
$G \times G$-stable subset $V$. But then $s_{\rm ad}$
vanishes on a closed $G \times G$-orbit in the 
closure of $V$ which can only be $\bm Y$. As above 
this is a contradiction and we conclude that 
$(s)_0 = c D_{i}$. 

In order to prove that $c=1$ we may assume that 
$G= G_{\rm sc}$ is simply connected and that 
$X = G$. In this case the statement is 
well known (cf. proof of 6.1.11 in \cite{BK}).
\end{proof}

\subsection{Sections of the dualizing sheaf}
\label{gs}
Let again $X$ be a smooth equivariant embedding of $G$.
For each $i=1, \dots, l$, there exists a unique 
$G_{\rm sc} \times G_{\rm sc}$-linearization of the 
line bundle $\mathcal{L}(D_i)$. The set of global 
sections of $\mathcal{L}(D_i)$ may then be regarded as 
a $G_{\rm sc} \times G_{\rm sc}$-module. We claim

\begin{prop}
\label{globalsections}
There exists a $G_{\rm sc} \times 
G_{\rm sc}$-equivariant morphism 
$$ \psi_i : {\rm H}(\omega_i)^* \boxtimes 
{\rm H}(\omega_i) \rightarrow
{\rm H}^0\big(X, \mathcal{L}(D_i)\big),$$
such that $\psi_i(u_{\omega_i}^+ \boxtimes 
v_{\omega_j}^+)$ is the canonical section
of $\mathcal{L}(D_i)$.
\end{prop}
\begin{proof}
When $X$ is toroidal this follows by Lemma 
\ref{sections}. For a general smooth 
embedding $X$ there exists by Zariski's 
main theorem (cf. proof of Prop.6.2.6 \cite{BK})
an open subset $X'$ of $X$ such that $X'$ is
a toroidal embedding of $G$ and such 
that the complement $X-X'$ has 
codimension $\geq 2$ in $X$. As the 
statement is invariant under replacing 
$X$ with $X'$ the result now follows.
\end{proof}

\subsubsection{The complete toroidal case}
\label{gs1} 

When $X$ is a complete toroidal embedding of $G$ 
we may even give more structure to the map 
$\psi_i$ given in Proposition \ref{globalsections}.
To this end, let $Y$  denote a closed  $G \times G$-orbit
in $X$ and consider the restriction map
$$ i_{|Y}^* : {\rm H}^0\big(X, \mathcal{L}(D_i)\big)
\rightarrow {\rm H}^0\big(Y, \mathcal{L}(D_i)_{|Y}\big).$$
By Lemma \ref{closedorbit2} and Lemma \ref{sections}
it follows that $\mathcal{L}(D_i)_{|Y} \simeq 
\mathcal{L}(\omega_i,-w_0 \omega_i)$. Consequently,
there exists a $G_{\rm sc} \times G_{\rm 
sc}$-equivariant isomorphism 
$$  {\rm H}^0\big(Y, \mathcal{L}(D_i)_{|Y}\big) \simeq 
{\rm H}(-w_0 \omega_i) \boxtimes {\rm H}(\omega_i).$$ 
By Lemma \ref{closedorbit2} the composition of $i_{|Y}^*$ 
with $\psi_i$ is nonzero and hence we obtain a
commutative $G_{\rm sc} \times G_{\rm 
sc}$-equivariant diagram 
\begin{equation}
\notag
\vcenter{
\xymatrix{
{\rm H}({\omega_{i}})^{*} \boxtimes {\rm H}(\omega_i) 
\ar[d]  \ar[r]^{\psi_{i}} & {\rm H}^0(X, \mathcal{L}(D_{i}) )
\ar[d]^{i_{|Y}^*}  & \\   
{\rm H}(-w_0 \omega_{i}) \boxtimes 
{\rm H}(\omega_{i}) \ar[r]^(.55)\simeq 
& {\rm H}(Y, \mathcal{L}(D_{i})_{|Y}) \\       
}}
\end{equation}
where the left vertical map is defined by some nonzero 
$G_{\rm sc}$-equivariant map $ {\rm H}({\omega_{i}})^{*} 
\rightarrow {\rm H}(-w_0 \omega_{i})$. Notice that 
by Frobenius reciprocity the map ${\rm H}
({\omega_{i}})^{*} \rightarrow {\rm H}(-w_0 \omega_{i})$
is defined uniquely up to a nonzero constant. 

\section{Frobenius Splitting}

In this section we collect a number of facts from
the theory of Frobenius splitting. The presentation
will be sketchy only stating the results which
we will need. For a more thorough, and closely 
related, presentation we refer to \cite{T}.

\subsection{Frobenius splitting}

  Let $X$ denote a scheme of finite type over an algebraically closed
  field $k$ of characteristic $p>0$. The \emph{absolute Frobenius
  morphism} on $X$ is the morphism $F : X \rightarrow X$ of schemes,
  which is the identity on the set of points and where the associated
  map of sheaves $$ F^{\sharp} : \mathcal{O}_{X}  \rightarrow F_{*}
  \mathcal{O}_{X}$$ is the $p$-th power map. We say that $X$ is
  \emph{Frobenius split} (or just F-split) if there exists a morphism $s \in
  \Hom_{\mathcal{O}_{X}} (F_{*} \mathcal{O}_{X}, \mathcal{O}_{X})$
  such that the composition $s \circ  F^{\sharp}$ is the identity map
  on $\mathcal{O}_{X}$.  

\subsection{Stable Frobenius splittings along divisors}
  Let $D$ denote an effective Cartier divisor on $X$ with associated
  line bundle $\mathcal{O}_{X}(D)$ and canonical section
  $\sigma_{D}$. We say that $X$ is \emph{stably Frobenius split along
    $D$} if there exists a positive integer $e$ and a morphism $$ s
  \in \Hom_{\mathcal{O}_{X}} (F_{*}^{e} \mathcal{O}_{X}(D),
  \mathcal{O}_{X}),$$ such that $s(\sigma_{D})=1$. In this case we say
  that $s$ is a stable Frobenius splitting of $X$ along $D$ of degree
  $e$. Notice that $X$ is Frobenius split exactly when there exists a
  stable Frobenius splitting of $X$ along the zero divisor $D=0$.

  \begin{remark}
    \label{remarkFsplittingonaopensubset}
    Consider an element $s \in \Hom_{\mathcal{O}_{X}} (F_{*}^{e}
    \mathcal{O}_{X}(D), \mathcal{O}_{X})$. Then the condition
    $s(\sigma_{D})=1$ on $s$ for it to define a stable Frobenius
    splitting of $X$, may be checked on any open dense subset of $X$.   
  \end{remark}

\subsection{Subdivisors}
\label{sectionsubdivisor}

  Let $D' \leq D$ denote an effective Cartier subdivisor and let $s$
  be a stable Frobenius splitting of $X$ along $D$ of degree $e$. The
  composition of $s$ with the map $$ F_{*}^{e} \mathcal{O}_{X} (D')
  \rightarrow F_{*}^{e} \mathcal{O}_{X}(D),$$ defined by the canonical
  section of the divisor $D-D'$, is then a stable Frobenius splitting
  of $X$ along $D'$ of degree $e$. Applying this to the case $D'=0$ it
  follows that if $X$ is stably Frobenius split along any effective
  divisor $D$ then $X$ is also Frobenius split.

\subsection{Compatibly split subschemes}
  Let $Y$ denote a closed subscheme of $X$ with sheaf of ideals
  $\mathcal{I}_{Y}$. When $$s \in \Hom_{\mathcal{O}_{X}} (F_{*}^{e}
  \mathcal{O}_{X}(D), \mathcal{O}_{X})$$ is a stable Frobenius
  splitting of $X$ along $D$ we say that $s$ \emph{compatibly
    Frobenius splits} $Y$ if the following conditions are satisfied 

  \begin{enumerate}
    \item The support of $D$ does not contain any of the irreducible
      components of $Y$. 
    \item $s\big(F_{*}^{e} (\mathcal{I}_{Y} \otimes \mathcal{O}_{X}(D))\big)
      \subseteq \mathcal{I}_{Y}$. 
  \end{enumerate}

When $s$ compatibly Frobenius splits $Y$
there exists an induced stable Frobenius 
splitting of $Y$ along $D \cap Y$ of degree
$e$.

  \begin{lemma}
    \label{5.4}
    Let $s$ denote a stable Frobenius splitting of $X$ along $D$ which
    compatibly Frobenius splits a closed subscheme $Y$ of $X$. If $D'
    \leq D$ then the induced stable Frobenius splitting of $X$ along
    $D'$, defined in Section~\ref{sectionsubdivisor}, compatibly
    Frobenius splits $Y$.  
  \end{lemma}

\begin{lemma}
\label{5.5}
Let $D_{1}$ and $D_{2}$ denote effective Cartier divisors. If
$s_{1}$ (resp. $s_{2}$) is a stable Frobenius splitting of $X$
along $D_{1}$ (resp. $D_{2}$) of degree $e_{1}$ (resp. $e_{2}$)
which compatibly splits a closed subscheme $Y$ of $X$, then there
exists a stable Frobenius splitting of $X$ along $D_{1} + D_{2}$
of degree $e_{1} + e_{2}$ which compatibly splits $Y$.
\end{lemma}

  \begin{lemma}
    \label{5.6}
    Let $s$ denote a stable Frobenius splitting of $X$ along an
    effective divisor $D$. Then 
    \begin{enumerate}
        \item If $s$ compatibly Frobenius splits a closed subscheme $Y$
        of $X$ then $Y$ is reduced and each irreducible component of 
        $Y$ is also
        compatibly Frobenius split by $s$. 
      \item Assume that $s$ compatibly Frobenius splits closed
        subschemes $Y_{1}$ and $Y_{2}$ and that the support of $D$
        does not contain any of the irreducible components of the
        scheme theoretic intersection $Y_{1} \cap Y_{2}$. Then $s$
        compatibly Frobenius splits  $Y_{1} \cap Y_{2}$.  
    \end{enumerate}
  \end{lemma}

  The following statement relates stable Frobenius splitting along
  divisors to compatibly Frobenius splitting.

  \begin{lemma}
    \label{5.7} 
    Let $D$ and $D'$ denote effective Cartier divisors and let $s$
    denote a stable Frobenius splitting of $X$ along $(p-1)D + D'$ of
    degree $1$. Then there exists a stable Frobenius splitting of $X$
    along $D'$ of degree $1$ which compatibly splits the closed
    subscheme defined by $D$.  
    
  \end{lemma}

\subsection{Cohomology and Frobenius splitting}

The notion of Frobenius splitting is particular useful in connection
with proving higher cohomology vanishing for line bundles. We will
need

\begin{lemma}
\label{lemmacohomologyandcompatiblyFsplit}
Let $s$ denote a stable Frobenius splitting of $X$ along $D$ of
degree $e$  and let $Y$ denote a closed compatibly Frobenius split
subscheme of $X$. Then for every line bundle $\mathcal L$ on $X$
and every integer $i$ there exists an inclusion 
$${\rm H}^{i}
(X, \mathcal{I}_{Y} \otimes \mathcal{L}) \subseteq {\rm H}^{i}
(X, \mathcal{I}_{Y} \otimes \mathcal{L}^{p^{e}} \otimes
\mathcal{O}_{X}(D)),$$ of abelian groups. In particular, when $X$
is projective, $\mathcal L$ is globally generated and $D$ is ample
then the group ${\rm H}^{i}(X, \mathcal{I}_{Y} \otimes
\mathcal{L})$ is zero for $i>0$. 
\end{lemma}
 
\subsection{Push forward}
\label{sectionpushforward}

Let $f : X \rightarrow X'$ denote a proper morphism of schemes and
assume that the induced map $ \mathcal{O}_{X'} \rightarrow f_{*}
\mathcal{O}_{X}$ is an isomorphism. Then every Frobenius splitting
of $X$ induces, by application of the functor $f_{*}$, a Frobenius
splitting of $X'$. Moreover, when $Y$ is a compatibly Frobenius
split subscheme of $X$ then the induced Frobenius splitting of $X'$
compatibly splits the scheme theoretic image $f(Y)$ (see 
\cite[Prop.4]{MR}). We will need the following connected statement.

\begin{lemma}
\label{lemmapushforward}
Let $f : X \rightarrow X'$ denote a morphism of projective schemes
such that $\mathcal{O}_{X'} \rightarrow f_{*} \mathcal{O}_{X}$ is
an isomorphism. Let $Y$ be a closed subscheme of $X$ and denote by
$Y'$ the scheme theoretic image $f(Y)$. Assume that there exists a
stable Frobenius splitting of $X$ along an ample divisor $D$ which
compatibly splits $Y$. Then $f_{*} \mathcal{O}_{Y} =
\mathcal{O}_{Y'}$ and ${\rm R}^{i} f_{*} \mathcal{O}_{Y} = 0$
for $i>0$. 
\end{lemma}

\subsection{Frobenius splitting of smooth varieties}

When $X$ is a smooth variety there exists a canonical $\mathcal{O}_{X}$-linear 
identification (see e.g. \cite[\S 1.3.7]{BK}) $$ F_{*}
\omega_{X}^{1-p} \simeq \Hom_{\mathcal{O}_{X}} (F_{*}
\mathcal{O}_{X}, \mathcal{O}_{X}).$$ Hence, a Frobenius splitting of
$X$ may be identified with a global section of $\omega_{X}^{1-p}$
with certain properties. A global section $\tau$ of $\omega_{X}^{1-p}$
which corresponds to a Frobenius splitting up to a nonzero 
constant will be called a \emph{Frobenius splitting section}. 

\begin{lemma}
\label{5.11}
Let $\tau$ be a Frobenius splitting section of a smooth 
variety $X$. Then there exists a stable Frobenius 
splitting of $X$ of
degree $1$ along the Cartier divisor defined by $\tau$. In
particular, if $\tau = \tilde{\tau}^{p-1}$ is a $(p-1)$-th power of
a global section $\tilde{\tau}$ of $\omega_{X}^{-1}$, then $X$ is
Frobenius split compatibly with the zero divisor of
$\tilde{\tau}$. 
\end{lemma}

\subsection{Frobenius splitting of $\nicefrac{\boldsymbol{G}}
{\boldsymbol{B}}$} 

The flag variety $X = \nicefrac{G}{B}$ is a smooth variety with 
dualizing sheaf  $\omega_X = \mathcal{L}(-2 \rho)$ where $\rho$ 
is a dominant weight defined as half of the sum of the positive 
roots. Let ${\rm St} := {\rm H}((p-1)\rho)$ denote the 
Steinberg module of $G_{\rm sc}$ and consider the multiplication
map 
$$ m_{\nicefrac{G}
{B}} : {\rm St} \otimes {\rm St} \rightarrow
  {\rm H}(2(p-1)\rho) \simeq 
 {\rm H}^0 (\nicefrac{G}{B}, \omega_{\nicefrac{G}{B}}^{1-p}).$$ 
The Steinberg module ${\rm St}$ is an irreducible
selfdual $G_{\rm sc}$-module and hence there exists 
a unique (up to nonzero scalars) nondegenerate 
$G_{\rm sc}$-invariant bilinear form 
$$ \phi_{{\nicefrac{G}{B}}} : {\rm St} \otimes {\rm St}
\rightarrow k.$$
Then (see \cite[Thm.2.3]{LT}) 

\begin{thm} 
\label{thmG/BFsplit}
Let $t$ be an element in ${\rm St} \otimes {\rm St}$. 
Then  $\phi_{\nicefrac{G}{B}} (t)$
is a Frobenius splitting section of $\nicefrac{G}{B}$ 
if and only if $\phi_{{\nicefrac{G}{B}}}(t)$ is nonzero.
\end{thm}

\section{F-splitting of smooth equivariant embeddings}
\label{F}
Let $X$ denote a smooth equivariant embedding of $G$.
Define $S$ to be the $G_{\rm sc} \times  G_{\rm sc}$-module 
$$S = \bigotimes _{i=1}^l \Big({\rm H}(\omega_i)^* 
\boxtimes{\rm H}(\omega_i)\Big)^{\otimes (p-1)}. $$
By Proposition \ref{globalsections} there exists a
$G_{\rm sc} \times  G_{\rm sc}$-equivariant morphism
$$\psi_X : S  \rightarrow  
{\rm H}^0\Big(X, {\mathcal L}\big( (p-1)\sum_{i=1}^l D_i\big) 
\Big). 
$$
defined as the $(p-1)$-th product of the $\psi_i$'s. 
Let $\sigma_j$ denote the canonical section of 
$\mathcal{L}(X_j)$, for $j=1, \dots, n$, and define 
for $s,t \in S$ the section
$$\Psi_X(s,t) = \psi_X(s) \psi_X(t) 
\prod_{i=1}^n \sigma_i^{p-1},$$
of the line bundle $\omega_X^{1-p}$ on $X$. 
Notice that if $X'$ is an equivariant embedding of 
$G$ which moreover is an open subset of $X$, then 
the restriction of $\Psi_X(s,t)$ to $X'$ is equal 
to $\Psi_{X'}(s,t)$.
The main result Theorem \ref{smoothembedding}
in this section describes 
when $\Psi_X(s,t)$ is a Frobenius splitting 
section of the smooth embedding $X$.

\subsection{F-splitting smooth complete toroidal embeddings} 
Consider a smooth complete toroidal embedding $X$ of $G$ and 
choose a closed $G \times G$-orbit $Y$ in $X$. By Lemma 
\ref{closedorbit2} we may identify $Y$ with $\nicefrac{G}{B} 
\times \nicefrac{G}{B}$ and under this isomorphism the restriction 
of $\mathcal{L}(D_i)$ to $Y$ corresponds to 
$\mathcal{L}(\omega_i, -w_0 \omega_i)$. In particular, restricting 
to $Y$ induces a map
$$ i_{| Y}^* : 
{\rm H}^0\Big(X, {\mathcal L}\big(2 (p-1)\sum_{i=1}^l D_i \big) \Big) 
\rightarrow {\rm H}^0\Big(Y, \omega_Y^{1-p} \Big).$$
This leads to the following result which also  
explains the standard way of Frobenius splitting $X$ 
(cf. proof of Thm.6.2.7 \cite{BK})

\begin{lemma}
\label{idea}
Let $X$ denote a smooth complete toroidal embedding of
$G$ and let 
$Y$ denote a closed $G \times G$-orbit in $X$. Let 
$s$ and $t$ be elements of $S$. Then $\Psi_X(s,t)$
is a Frobenius splitting section of $X$ if and 
only if the restriction of $\psi_X(s) \psi_X(t)$
to $Y$ is a Frobenius splitting section of $Y$.
\end{lemma}

In order to control the restriction of  $\psi_X(s) \psi_X(t)$
to $Y$ we use Section \ref{gs1}. It follows that we have a 
commutative $G_{\rm sc} \times G_{\rm sc}$-equivariant 
diagram with nonzero maps  
$$
\xymatrix 
{S \ar[r]^(.25){\psi_X} 
\ar[d]_\eta &
{\rm H}^0\Big(X, {\mathcal L}\big( (p-1)\sum_{i=1}^l D_i\big) 
\Big) \ar[d]^{i_{|Y}^*}  
 \\
{\rm St} \boxtimes{\rm St} \ar[r]^(.45){\simeq} & {\rm H}^0\Big(Y, \omega_Y^{
\nicefrac{(1-p)}{2}} \Big)
}
$$
where ${\rm St} = {\rm H}\big((p-1) \rho\big)$ denotes 
the Steinberg module of $G_{\rm sc}$. 
Using the $G_{\rm sc} \times G_{\rm sc}$-invariant 
form $\phi_{\nicefrac{G}{B} \times \nicefrac{G}{B}}$ 
on ${\rm St} \boxtimes {\rm St}$ we may define a similar 
form on $S$  by 
$$ \phi : S \otimes S \rightarrow k,$$
$$ s \otimes t \mapsto 
\phi_{\nicefrac{G}{B} \times \nicefrac{G}{B}}(\eta(s) \otimes \eta(t))$$
Notice that $S$ and the  $G_{\rm sc} \times 
G_{\rm sc}$-invariant form $\phi$ is defined 
without the help of $X$. In particular, 
$S$ and $\phi$ does not depend on $X$.
Now by Lemma \ref{idea} and Theorem
\ref{thmG/BFsplit} we find 

\begin{prop}
\label{completetoroidal}
Let the notation be as above and let $s$ and $t$ be 
elements of $S$. Then 
$\Psi_X(s,t)$
is a Frobenius splitting section of $X$ if and 
only if $\phi(s \otimes t)$ is nonzero.   
\end{prop}

\subsection{Frobenius splitting $G$} 

By restricting the statement of Proposition 
\ref{completetoroidal} to $G$ we find 

\begin{cor}
\label{group}
Let $s$ and $t$ be elements of $S$. Then 
$$\Psi_G(s,t) :=\psi_G(s) \psi_G(t), $$
is a Frobenius splitting section of $G$ if 
and only if $\phi(s \otimes t)$ is nonzero.   
\end{cor}
\begin{proof}
Choose a smooth complete toroidal embedding $X$ 
of $G$ and consider $\Psi_X(s,t)$. Remember 
that checking whether $\Psi_X(s,t)$ is a 
Frobenius splitting section of $X$ may be 
done on the open subset $G$ (see Remark 
\ref{remarkFsplittingonaopensubset}). 
Now apply Proposition \ref{completetoroidal}.
\end{proof}

\subsection{Frobenius splitting smooth 
equivariant embeddings}
\label{see}

We can now prove that main result of this
section.

\begin{thm}
\label{smoothembedding}
Let $X$ denote an arbitrary smooth embedding 
of $G$ and let $s$ and $t$ be elements of $S$. 
Then $\Psi_X(s,t)$
is a Frobenius splitting section of $X$ if 
and only if $\phi(s \otimes t)$ is nonzero.   
\end{thm}
\begin{proof}
That $\Psi_X(s,t)$ is a Frobenius splitting
section may be checked on the open subset 
$G$. Now apply Corollary \ref{group}.
\end{proof}

In the following statement $t_i$, $i=1,\dots,l$, denotes
the identity map in ${\rm End}({\rm H}(\omega_i)^*) \simeq 
{\rm H}(\omega_i)^* \boxtimes{\rm H}(\omega_i)$. Notice
that as an element of ${\rm End}({\rm H}(\omega_i))^*$ 
the element $t_i$ is just the trace map on 
${\rm End}({\rm H}(\omega_i))$. We also fix a 
nonzero weight vector $u_{\omega_i}^-$
of ${\rm H}(\omega_i)^*$ of weight $-\omega_i$. 

\begin{cor}
\label{corsmooth}
The global section 
$$ \prod_{i=1}^l \psi_i(t_i)^{p-1} \prod_{i=1}^l \psi_i
(u_{\omega_i}^-  \boxtimes v_{\omega_i}^+)^{p-1} \prod_{j=1}^n
\sigma_j^{p-1},$$
of $\omega_x^{1-p}$ is a Frobenius splitting section 
of $X$. 
\end{cor}
\begin{proof}
It suffices by Theorem \ref{smoothembedding} to prove that 
$$\phi\Big(\bigotimes_{i=1}^l t_i^{\otimes (p-1)} \otimes \bigotimes_{i=1}^l 
(u_{\omega_i}^-  \boxtimes v_{\omega_i}^+)^{\otimes (p-1)}\Big) $$
is nonzero. The image of $\bigotimes_{i=1}^l t_i^{\otimes (p-1)}$
in ${\rm St} \boxtimes {\rm St}$ 
coincides with a nonzero diag$(G)$-invariant 
element $v_\Delta$. Moreover, the image of the element
$ \bigotimes_{i=1}^l (u_{\omega_i}^-  \boxtimes 
v_{\omega_i}^+)^{\otimes (p-1)} $ in ${\rm St} \boxtimes {\rm St}$ 
equals $v_- \boxtimes v_+$ for some nonzero weight 
vectors $v_+$ and $v_-$ in ${\rm St}$ of weight
$(p-1) \rho$ and $-(p-1)\rho$ respectively. Thus,
we have to show that $\phi_{\nicefrac{G}{B} \times \nicefrac{G}{B}}
\big(v_\Delta \otimes  (v_- \boxtimes v_+)\big)$
is nonzero. But by weight consideration this is clearly the 
case. 
\end{proof}

\section{Consequences in the smooth case}

In this section we collect a number of consequences of the 
results in Section \ref{F} and the following Lemma 
\ref{lemmaFsplitgivespurecodim}.
Notice that when $f$ is a global section of a line bundle
$\mathcal{L}$ on a variety $X$, then we may regard $f$ 
as an element in the local rings $\mathcal{O}_{X,x}$ at 
points $x \in X$. This identification is unique up to 
units in $\mathcal{O}_{X,x}$. Using this identification
we may now state

\begin{lemma}
\label{lemmaFsplitgivespurecodim}
Let $X$ denote a smooth variety with dualizing sheaf
$\omega_{X}$ and let  $\mathcal{L}_{1}, \dots , 
\mathcal{L}_{N}$ denote a collection of line bundles
on $X$ such that $\otimes_{i=1}^{N} \mathcal{L}_{i}
\simeq \omega_{X}^{-1}$. Let $f_i$, $i=1,\dots,N$, 
denote a global section of $\mathcal{L}_i$ and 
assume that $\prod_{i=1}^{N} f_{i}^{p-1}$, 
considered as a global section of $\omega_X^{1-p}$, 
is a Frobenius splitting section of $X$. Choose
a sequence $1 \leq i_1, \dots, i_r \leq N$ of pairwise
distinct integers. Then  
\begin{enumerate}
\item  The sequence $f_{i_{1}}, \dots , 
f_{i_{r}}$ forms a regular sequence in the local ring
$\mathcal{O}_{X,x}$ at a point $x$ contained in the 
common zero set of $f_{i_1}, \dots, f_{i_r}$. 
\item  
The common zero set of $f_{i_1}, \dots, f_{i_r}$ has 
pure codimension $r$.
\end{enumerate}
\end{lemma}
\begin{proof}
As all statements are local we may assume that $X$
is affine and that $\omega_X$ and $\mathcal{L}_1,
\dots, \mathcal{L}_N$ are all trivial. Hence, the
elements $f_1, \dots, f_N$ are just regular global 
functions on $X$. Moreover, by assumption there exists
a function $\tau : F_* k[X] \rightarrow k[X]$ such 
that $\tau\big(a^p (f_1 \cdots f_N)^{p-1}\big) = a$ for 
all global regular functions $a \in k[X]$.

Let $x$ be a common zero of 
$f_{i_1}, \dots, f_{i_r}$ and assume that we
have a relation of the form $\sum_{s=1}^{j} a_{s} f_{i_{s}} 
= 0$ for certain elements $a_s$ in  $\mathcal{O}_{X,x}$. 
In particular, the product 
$a_{j}^{p} (f_1 \cdots f_N)^{p-1}$
is contained in the ideal $(f_{i_1}^p, \dots, f_{i_{j-1}}^p)$
of $\mathcal{O}_{X,x}$ and hence 
$$ a_j = \tau\big(a_{j}^{p} (f_1 \cdots f_N)^{p-1}\big) 
\in (f_{i_1}, \dots, f_{i_{j-1}}).$$
This proves (1). Now (2) follows as a direct consequence 
of (1).
\end{proof}

We can now prove the first of our main results

\begin{cor}
\label{geometryN}
Let $X$ be a smooth equivariant embedding of $G$
and let $\bar{N}$ denote the closure of the
Steinberg zero-fiber in $X$. Then 
\begin{enumerate}
\item $\bar{N}$ coincides with the scheme theoretic 
intersection of the zero sets of $t_i$, $i=1, \dots, l$.
In particular, $\bar{N}$ is a local complete intersection.
\item $\bar{N}$ is normal, Gorenstein and Cohen-Macaulay.
\item The dualizing sheaf of $\bar{N}$ is isomorphic
to the restriction of the line bundle 
$\mathcal{L}_{\bar{N}} := \mathcal{L}(-\sum_{i=1}^l (w_0,1) D_i 
- \sum_{j=1}^n X_j)$
to $\bar{N}$.   
\end{enumerate}
\end{cor}
\begin{proof}
Consider the Frobenius splitting section of $X$ 
defined in Corollary \ref{corsmooth}. By 
Lemma \ref{5.11} and Lemma \ref{5.6} the scheme theoretic intersection
$C$ of  $t_i$, $i=1, \dots, l$, is a reduced
scheme. Moreover, by  Lemma \ref{lemmaFsplitgivespurecodim} 
each component of $C$ has codimension $l$
and will intersect the open locus $G$ (else,
by Lemma  \ref{lemmaFsplitgivespurecodim}, 
such a component would have codimension $\geq l+1$). 
We conclude that $C \cap G$ is dense in $C$
and that $C$ is a local complete intersection.
But clearly (see remark above Corollary 
\ref{corsmooth}) $C \cap G$ coincides with
the Steinberg zero-fiber $N$, and thus $C$ must
be equal to the closure $\bar{N}$. This 
proves (1). 

To prove (2) it then suffices to show that 
$\bar{N}$ is regular in codimension 1. 
Let $Z$ denote a component of the singular 
locus of $\bar{N}$. If $Z \cap G \neq 
\emptyset$ then the codimension of $Z$ is 
$\geq2$ as $N$ is normal by \cite[Thm.6.11]{St}. 
So assume that $Z$ is contained 
in a boundary component $X_{j}$. Now, 
by Lemma \ref{5.11} the scheme theoretic 
intersection $\bar{N} \cap X_{j}$ is reduced. 
Hence, as $X_j$ is a Cartier divisor, every 
smooth point of $\bar{N} \cap X_{j}$ is also 
a smooth point of $\bar{N}$. In particular, 
$Z$ is properly contained in a component of
$\bar{N} \cap X_{j}$. But the variety
$\bar{N} \cap X_{j}$ has pure codimension 1 
in $\bar{N}$ which ends the proof of (2).

Statement (3) follows by (1) and the  
description of the dualizing sheaf of $X$ 
in Proposition \ref{canonical}.
\end{proof}

\subsection{Stable Frobenius splittings along divisors}

\begin{prop}
\label{propNbar,stableFsplitting}
Let $X$ be a smooth equivariant embedding of $G$. Then there exists a
stable Frobenius splitting of $X$ along the divisor
\begin{displaymath}
(p-1) \Big( \sum_{j=1}^{n} X_{j} + \sum_{i=1}^{\ell} (w_0,1)D_{i} \Big)
\end{displaymath}
of degree 1 which compatibly Frobenius splits the closure 
$\bar{N}$ of the Steinberg zero-fiber. 
\end{prop}
\begin{proof}
Let $\tau$ denote the Frobenius splitting section of 
Corollary \ref{corsmooth}. By Lemma \ref{5.11},
Lemma \ref{5.7} and Lemma \ref{globalsections} we know that $\tau$ defines 
a degree 1 stable Frobenius splitting of $X$ along the divisor 
$$
(p-1) \Big( \sum_{j=1}^{n} X_{j} + \sum_{i=1}^{\ell} (w_0,1)D_{i} \Big),  
$$
which compatibly Frobenius splits the zero divisor of 
$ \prod_{i=1}^l \psi_i(t_i)$. Now apply Lemma 
\ref{5.6}(2), Corollary \ref{geometryN} and 
Lemma \ref{lemmaFsplitgivespurecodim}. 
\end{proof}

\begin{cor}
\label{corNbarstableFsplitalongampledivisor}
Let $X$ denote a projective smooth equivariant embedding
of $G$. Then there exists a stable Frobenius splitting of 
$X$ along an ample divisor with support $X - G$ 
which compatibly Frobenius splits the subvariety $\bar{N}$.
\end{cor}
\begin{proof}
By Proposition \ref{propNbar,stableFsplitting} and
Lemma \ref{5.4} there exists a stable
Frobenius splitting of $X$ along the divisor $\sum_{j=1}^{n}
X_{j}$ which compatibly splits $\bar{N}$. Applying
Lemma \ref{5.4} and
Lemma \ref{5.5} it suffices to
show that there exist positive integers $c_{j}>0$ such that
$\sum_{j=1}^{n} c_{j} X_{j}$ is ample. This follows from
\cite[Prop.4.1(2)]{BT}. 
\end{proof}

This has the following implications for resolutions 

\begin{cor}
\label{cordirectimages}
Let $X$ be a projective equivariant embedding of $G$
and let $f: X' \to X$ be a projective resolution of 
$X$ by a smooth projective equivariant $G$-embedding $X'$. 
Denote by $\bar{N}'$ (resp. $\bar{N}$) the closure of the 
Steinberg zero-fiber within $X'$ (resp. $X$). Then \\
(i) $f_{*} \mathcal{O}_{X'} = \mathcal{O}_{X}$ and ${\rm R}^{i}
f_{*} \mathcal{O}_{X'} = 0$ for $i > 0$. (cf. \cite[pf. of Cor.2]{Rit})\\
(ii) $f_{*} \mathcal{O}_{\bar{N}'} = \mathcal{O}_{\bar{N}}$ and
${\rm R}^{i} f_{*} \mathcal{O}_{\bar{N}'} = 0$ for $i > 0$. 

\end{cor}
\begin{proof}
As $X'$ is normal and $f$ is birational it follows from Zariski's
main theorem that $f_{*} \mathcal{O}_{X} =
\mathcal{O}_{X'}$. Hence, by Lemma~\ref{lemmapushforward} it
suffices to prove that there exists a stable Frobenius splitting
of $X$ along an ample divisor which compatibly Frobenius splits
$\bar{N}$. Now apply
Corollary~\ref{corNbarstableFsplitalongampledivisor}. 
\end{proof}

\section{Frobenius splitting $\bar{N}$ for general embeddings}

In this section $X$ will denote an arbitrary equivariant embedding
of $G$ and $\bar{N}$ will denote the closure of the Steinberg 
zero-fiber in $X$. 

\begin{thm}
\label{thmNbarFsplit}
There exists a Frobenius splitting of $X$ which simultaneously 
compatibly splits the closed subvarieties
$\bar{N}, (w_0,1){D}_{i}, X_{j}$, for $i=1, \dots , l$ and
$j=1, \dots, n$. 
\end{thm}
\begin{proof}
By Theorem \ref{thmprojectiveresolution} we may 
find a projective resolution $f : X' \to X$ by 
a smooth toroidal embedding  $X'$ of $G$. 
By Zariski's Main Theorem, $f_{*} \mathcal{O}_{X'} =
\mathcal{O}_{X}$. Thus, by \cite[Prop.4]{MR} 
(cf. section~\ref{sectionpushforward}) we can reduce to the case
where $X$ is smooth. Now apply Corollary \ref{corsmooth},
Lemma \ref{5.11}, Lemma \ref{5.6}, Lemma \ref{globalsections}
and Corollary \ref{geometryN} in the given order.
\end{proof}

\begin{example}
\label{example}
Consider the group $G = {\rm PSL}_2(k)$ over a field
$k$ of positive characteristic different from $2$. 
Then the wonderful compactification ${\bm X}$ 
of $G$ may be identified with the projectivization 
of the set of $2 \times 2$-matrices with entries in $k$. 
Denote the homogeneous coordinates in $\bm X$
by $a,b,c$ and $d$. Then the closure $\bar{\mathcal{U}}$ 
of the unipotent variety $\mathcal{U}$ of $G$
within $\bm X$, is defined by the polynomial 
$f=(a+d)^2 - 4 (ad-bc)$. Moreover, the boundary is 
defined by the polynomial $g=(ad-bc)$.
In particular, the ideal generated by $f$ and $g$
is not reduced and, as a consequence, the boundary 
${\bm X}-G$ and the closure $\bar{\mathcal{U}}$ 
cannot be compatibly Frobenius split at the same
time. When $k$ has characteristic $2$ the unipotent
variety $\mathcal{U}$ coincides with the Steinberg
zero-fiber. In this case the polynomial defining
$\bar{\mathcal{U}}$ is given by $f=a+d$ and we 
do not see a similar problem. 
\end{example}

\begin{remark}
W. van der Kallen and T. Springer has informed us
that they have proved Theorem \ref{thmNbarFsplit} 
in case $X$ is the wonderful compactification of 
a group of adjoint type. Their proof proceeds by
descending the Frobenius splitting results in 
\cite{T} to the wonderful compactification. 
\end{remark}

We can also prove a vanishing result for line
bundles on $\bar{N}$ :

\begin{prop}
\label{propcohomologyonNbar}
Let $X$ denote a projective equivariant $G$-embedding and let
$\mathcal{L}$ (resp. $\mathcal{M}$) denote a globally generated 
line bundle on $X$ (resp. $\bar{N}$). Then 
$${\rm H}^i(X, \mathcal{L}) =
{\rm H}^i(\bar{N}, \mathcal{M}) = 0, \, i>0.$$ 
Moreover, the restriction
map $$\Ho{X, \mathcal{L}} \to \Ho{\bar{N}, \mathcal{L}},$$ is
surjective. 
\begin{proof}
By Corollary~\ref{cordirectimages} we may assume that $X$ is
smooth. Now apply
Corollary~\ref{corNbarstableFsplitalongampledivisor} and the ``in
particular'' part of Lemma~\ref{lemmacohomologyandcompatiblyFsplit}.  
\end{proof}
\end{prop}

\subsection{Canonical Frobenius splittings of $X$}
A Frobenius splitting $s: F_{*} \mathcal{O}_{Z} \to
\mathcal{O}_{Z}$ of a $B$-variety $Z$ is a $T$-invariant
Frobenius splitting such that the action of a root 
subgroup of $G$ associated to the simple root $\alpha_i$, 
is of the form
$$x_{\alpha_{i}}(c) s = \sum_{j=1}^{p-1} c^{j} s_{j},$$ 
for certain morphisms $s_{j}: F_{*}
\mathcal{O}_{Z} \to \mathcal{O}_{Z}$ and all $c \in k$.

As a subset of $X$ the closure $\bar{N}$ is invariant 
under the diagonal action of $G$. In particular,
$\bar{N}$ is invariant under diag$(B)$ and we claim

\begin{lemma}
The variety $\bar{N}$ is canonical Frobenius split
with respect to the action of ${\rm diag}(B)$.
\end{lemma}
\begin{proof}
It suffices to prove that $X$ has a diag$(B)$-canonical 
splitting which compatibly splits $\bar{N}$. Moreover, 
by Theorem \ref{thmprojectiveresolution} we may assume
that $X$ is smooth. By the proof of Corollary 
\ref{geometryN} it then suffices to prove that the
Frobenius splitting section of Corollary \ref{corsmooth}
is canonical.  

As $\psi_{i}^{*} (t_{i})$ and $\sigma_{j}$ are
diag$(G)$-invariant we may concentrate on the
diag$(T)$-invariant factors $\psi_{i}^{*} (u_{\omega_i}^{-} 
\otimes v_{\omega_i}^{+})$. The statement follows now as 
\begin{equation*}
\begin{alignedat}{1}
& x_{\alpha_{j}}(c). v_{\omega_i}^{+} = v_{\omega_i}^{+}, \\
& x_{\alpha_{j}}(c). u_{\omega_i}^{-} =  u_{\omega_i}^{-}+ c u_{i,j} \ ,
\end{alignedat}
\end{equation*}
for certain elements $u_{i,j} \in {\rm H}(-w_{o} \omega_{i})^{*}$. 
\end{proof}
 
As a consequence we have 
(see \cite[Thm.4.2.13]{BK})

\begin{prop}
Let $\mathcal{L}$ denote a $G_{\rm sc}$-linearized 
line bundle on $\bar{N}$. Then the $G_{\rm sc}$-module
${\rm H}^0(\bar{N}, \mathcal{L})$ admits a good filtration,
i.e. there exists a filtration by $G_{\rm sc}$-modules 
$$0 =M^0 \subseteq M^1 \subseteq M^2 \subseteq \cdots 
\subseteq {\rm H}^0(\bar{N}, \mathcal{L}),$$ 
such that ${\rm H}^0(\bar{N}, \mathcal{L}) = 
\cup_i M^i$ and satisfying that the successive quotients
$M^{j+1}/M^{j}$ are isomorphic to modules of the form
${\rm H}(\lambda_j)$ for certain dominant weights $\lambda_j$. 
\end{prop}

\section{Geometric properties of $\bar{N}$}

Let $X$ be an arbitrary equivariant $G$-embedding. When 
$X$ is smooth we have seen that $\bar{N}$ 
is normal and Cohen-Macaulay. 
In this section we will extend these two properties to
arbitrary equivariant embeddings. 

The following result is due to G. Kempf although the version 
stated here is taken from \cite[\S 7]{BP} :

\begin{lemma}
\label{lemmaKempf}
Let $f: Z' \to Z$ denote a proper map of algebraic schemes
satisfying that $f_{*} \mathcal{O}_{Z'} \simeq \mathcal{O}_{Z}$
and ${\rm R}^{i} f_{*} \mathcal{O}_{Z'} = 0$ for $i > 0$. If
$Z'$ is Cohen-Macaulay with dualizing sheaf $\omega_{Z'}$ and
if ${\rm R}^{i} f_{*} \omega_{Z'} = 0$ for $i>0$, then $Z$ is
Cohen-Macaulay with dualizing sheaf $f_{*} \omega_{Z'}$.  
\end{lemma}

We will also need the following result due to V. Mehta and W. van
der Kallen (\cite[Thm.1.1]{MvdK}):

\begin{lemma}
\label{lemmaMvdK}
Let $f: Z' \to Z$ denote a proper morphism of schemes and let $V'$
(resp. $V$) denote a closed subscheme of $Z'$ (resp. $Z$). By
$\mathcal{I}_{V'}$ we denote the sheaf of ideals of $V'$. Fix an
integer $i$ and assume
\begin{enumerate}
\item $f^{-1}(V) \subseteq V'$.
\item ${\rm R}^{i} f_{*} \mathcal{I}_{V'}$ vanishes outside $V$.
\item $V'$ is compatibly F-split in $Z'$.
\end{enumerate}
Then ${\rm R}^{i} f_{*} \mathcal{I}_{V'} =0$.
\end{lemma}

We are ready to prove

\begin{thm}
Let $X$ denote an arbitrary equivariant $G$-embedding. Then the
closure $\bar{N}$ of the Steinberg zero-fiber in $X$ is normal and
Cohen-Macaulay. 
\end{thm}
\begin{proof}
Any equivariant embedding has an open cover by open equivariant
subsets of projective equivariant embeddings (see e.g. proof of
\cite{BK} Corollary 6.2.8). This reduces the statement to the
case where $X$ is projective. Choose a projective resolution $f
: X' \rightarrow X$ of $X$ by a smooth equivariant embedding
$X'$. By Corollary \ref{cordirectimages} we know that 
$f_{*} \mathcal{O}_{\bar{N}'} = \mathcal{O}_{\bar{N}}$
and applying Corollary~\ref{geometryN} this implies
that $\bar{N}$ is normal.

In order to show that $\bar{N}$ is Cohen-Macaulay we apply the
above Lemma \ref{lemmaKempf} and Lemma \ref{lemmaMvdK}. By
Corollary \ref{cordirectimages} it suffices to prove that 
${\rm R}^{i} f_{*} \omega_{\bar{N}'}= 0$, $i>0$, where 
$\omega_{\bar{N}'}$ is the dualizing sheaf of $\bar{N}'$.
By Corollary \ref{geometryN} the dualizing sheaf 
$\omega_{\bar{N}'}$ is isomorphic to the restriction
of $\mathcal{L}_{\bar{N}'}$ to $\bar{N}'$. Let 
$s'$ denote the canonical section of the line bundle
$\mathcal{L}_{\bar{N}'}^{-1}$ on $X'$ and let $V'$ denote 
the intersection of $\bar{N}'$ with the zero 
divisor of $s'$. Combining
Proposition \ref{propNbar,stableFsplitting} and
Lemma \ref{5.7} we find
that $\bar{N}'$ is Frobenius split compatibly with the closed
subscheme $V'$. Moreover, $f : \bar{N}' \to \bar{N}$ is an
isomorphism above the open subset $N$ and $$f^{-1}(\bar{N}
\setminus N) \subseteq V'.$$ Hence, by Lemma~\ref{lemmaMvdK} we
conclude ${\rm R}^{i} f_{*} \mathcal{I}_{V'}= 0$ for
$i>0$. But  $\mathcal{I}_{V'}$ is isomorphic to the 
restriction of $\mathcal{L}_{\bar{N}'}$ to
$\bar{N}'$. 
\end{proof}

\end{document}